\documentclass[12pt,reqno]{article}

\usepackage{amsmath,amsthm,amsfonts,amssymb,latexsym,mathrsfs}

\newtheorem{thm}{Theorem}[section]

\newtheorem{cor}[thm]{Corollary}
\newtheorem{prop}[thm]{Proposition}
\newtheorem{conj}[thm]{Conjecture}

\newtheorem*{DP}{Davenport-P\'olya Theorem}

\theoremstyle{definition}

\newtheorem{ex}[thm]{Example}
\newtheorem{rem}[thm]{Remark}

\numberwithin{equation}{section}

\newcommand{\lrf}[1]{\lfloor #1\rfloor}

\newcommand{\la}{\lambda}
\newcommand{\A}{\mathscr{A}}

\title{On the log-convexity of combinatorial sequences
\thanks{Partially supported by NSF of China 10471016.}}
\author{Lily L. Liu,\quad Yi Wang
\footnote{Corresponding author.
\newline\hspace*{5mm}
   {\it Email addresses:}\quad lliulily@yahoo.com.cn (L.L. Liu), wangyi@dlut.edu.cn (Y. Wang).}}
\date{\footnotesize Department of Applied Mathematics,
         Dalian University of Technology,
         Dalian 116024, PR China}

\begin{document}

\maketitle

\begin{abstract}
This paper is devoted to the study of the log-convexity of
combinatorial sequences. We show that the log-convexity is preserved
under componentwise sum, under binomial convolution, and by the
linear transformations given by the matrices of binomial
coefficients and Stirling numbers of two kinds. We develop
techniques for dealing with the log-convexity of sequences
satisfying a three-term recurrence. We also introduce the concept of
$q$-log-convexity and establish the connection with linear
transformations preserving the log-convexity. As applications of our
results, we prove the log-convexity and $q$-log-convexity of many
famous combinatorial sequences of numbers and polynomials.
\end{abstract}

\section{Introduction}
\hspace*{\parindent}
Let $a_0,a_1,a_2,\ldots$ be a sequence of nonnegative numbers. The
sequence is called {\it convex} (resp. {\it concave}) if for $k\ge
1$, $a_{k-1}+a_{k+1}\ge 2a_k$ (resp. $a_{k-1}+a_{k+1}\le 2a_k$). The
sequence is called {\it log-convex} (resp. {\it log-concave}) if for
all $k\ge 1$, $a_{k-1}a_{k+1}\ge a_k^2$ (resp. $a_{k-1}a_{k+1}\le
a_k^2$). By the arithmetic-geometric mean inequality, the
log-convexity implies the convexity and the concavity implies the
log-concavity. Clearly, a sequence $\{a_k\}_{k\ge 0}$ of positive
numbers is log-convex (resp. log-concave) if and only if the
sequence $\left\{a_{k+1}/a_k\right\}_{k\ge 0}$ is increasing (resp.
decreasing).

It is well known that the binomial coefficients $\binom{n}{k}$, the
Eulerian numbers $A(n,k)$, the Stirling numbers $c(n,k)$ and
$S(n,k)$ of two kinds are log-concave in $k$ for fixed $n$
respectively (see \cite{WYjcta05} for instance). In contrast, it is
not so well known that many famous sequences in combinatorics,
including the Bell numbers, the Catalan numbers and the Motzkin
numbers, are log-convex respectively. Although the log-convexity of
a sequence of positive numbers is equivalent to the log-concavity of
its reciprocal sequence, they are also fundamentally different. For
example, the log-concavity of sequences is preserved by both
ordinary and binomial convolutions (see Wang and Yeh~\cite{WYjcta07}
for instance). But the ordinary convolution of two log-convex
sequences need not be log-convex. Even the sequence of partial sums
of a log-convex sequence is not log-convex in general. On the other
hand, Davenport and P\'olya~\cite{DP49} showed that the
log-convexity is preserved under the binomial convolution.

There have been quite a few papers concerned with the log-concavity
of sequences (see the survey articles \cite{Sta89,Bre94} and some
recent developments
\cite{LW-RZP,wyjcta02,wyeujc02,wylaa03,WYeujc05,WYjcta05,WYjcta07}).
However, there is no systematic study of the log-convexity of
sequences. Log-convexity is, in a sense, more challenging property
than log-concavity. One possible reason for this is that the
log-concavity of sequences is implied by the P\'olya frequency
property. An infinite sequence $a_0,a_1,a_2,\ldots$ is called a {\it
P\'olya frequency} (PF, for short) sequence if all minors of the
infinite Toeplitz matrix $(a_{i-j})_{i,j\ge 0}$ are nonnegative,
where $a_k=0$ if $k<0$. A finite sequence
$a_0,a_1,\ldots,a_{n-1},a_n$ is PF if the infinite sequence
$a_0,a_1,\ldots,a_{n-1},a_n,0,0,\ldots$ is PF. Clearly, a PF
sequence is log-concave. PF sequences are much better behaved and
have been studied deeply in the theory of total positivity
(Karlin~\cite{Kar68}). For example, the fundamental representation
theorem of Schoenberg and Edrei states that a sequence
$a_0=1,a_1,a_2,\ldots$ of real numbers is PF if and only if its
generating function has the form
$$\sum_{n\ge 0}a_nz^n
=\frac{\prod_{j\ge 1}(1+\alpha_jz)}{\prod_{j\ge
1}(1-\beta_jz)}e^{\gamma z}$$ in some open disk centered at the
origin, where $\alpha_j,\beta_j,\gamma\ge 0$ and $\sum_{j\ge
1}(\alpha_j+\beta_j)<+\infty$ (see Karlin~\cite[p.~412]{Kar68} for
instance). In particular, a finite sequence of nonnegative numbers
is PF if and only if its generating function has only real zeros
(\cite[p.~399]{Kar68}). So it is often more convenient to show that
a sequence is PF even if our original interest is only in the
log-concavity. Indeed, many log-concave sequences arising in
combinatorics are actually PF sequences.
Brenti~\cite{Bre89,Bre95,Bre96} has successfully applied total
positivity techniques and results to study the log-concavity
problems.

This paper is devoted to the study of the log-convexity of
combinatorial sequences and is organized as follows. In \S2 we
consider various operators on sequences that preserve the
log-convexity. We show that log-convexity is preserved under
componentwise sum, under binomial convolution, and by the linear
transformations given by the matrices of binomial coefficients,
Stirling numbers of the second kind, and signless Stirling numbers
of the first kind. In \S3 we discuss the log-convexity of sequences
satisfying a three-term recurrence. As consequences, some famous
combinatorial sequences, including the central binomial
coefficients, the Catalan numbers, the Motzkin numbers, the Fine
numbers, the central Delannoy numbers, the little and large
Schr\"oder numbers, are log-convex. And also, each of bisections of
combinatorial sequences satisfying three-term linear recurrences,
including the Fibonacci, Lucas and Pell numbers, must be log-convex
or log-concave respectively. In \S4 we introduce the concept of the
$q$-log-convexity of sequences of polynomials in $q$ and show the
$q$-log-convexity of the Bell polynomials, the Eulerian polynomials,
the $q$-Schr\"oder numbers and the $q$-central Delannoy numbers. We
also present certain linear transformations preserving the
log-convexity of sequences and establish the connection with the
$q$-log-convexity. Finally, in \S5, we present some conjectures and
open problems.
\section{Operators preserving log-convexity}
\hspace*{\parindent}
In this section we consider operators on sequences that preserve the
log-convexity. A similar problem for the log-concavity has been
studied (see \cite{Bre89,wylaa03,WYjcta07} for instance). However,
there are fundamentally different between them. For example, it is
somewhat surprising that the log-convexity is preserved under
componentwise sum.
\begin{prop}\label{lcx+lcx}
If both $\{x_n\}$ and $\{y_n\}$ are log-convex, then so is the
sequence $\{x_n+y_n\}$.
\end{prop}
\begin{proof} It follows immediately that
$$(x_{n-1}+y_{n-1})(x_{n+1}+y_{n+1})
\ge (\sqrt{x_{n-1}x_{n+1}}+\sqrt{y_{n-1}y_{n+1}})^2 \ge
(x_n+y_n)^2$$ from the well-known Cauchy's inequality and the
log-convexity of $\{x_n\}$ and $\{y_n\}$. \end{proof}

Given two sequences $\{x_n\}_{n\ge 0}$ and $\{y_n\}_{n\ge 0}$,
define their ordinary convolution by
$$z_n=\sum_{k=0}^{n}x_ky_{n-k}$$
and binomial convolution by
$$z_n=\sum_{k=0}^{n}\binom{n}{k}x_ky_{n-k},$$
respectively. It is known that the log-concavity of sequences is
preserved by both ordinary and binomial convolutions (see Wang and
Yeh~\cite{WYjcta07} for instance). However, the ordinary convolution
of two log-convex sequences need not be log-convex. On the other
hand, the binomial convolution of two log-convex sequences is
log-convex by Davenport and  P\'olya~\cite{DP49}.
Proposition~\ref{lcx+lcx} can provide an interpretation of this
result.
\begin{DP}\label{dp-bc}
If both $\{x_n\}$ and $\{y_n\}$ are log-convex, then so is their
binomial convolution
$$z_n=\sum_{k=0}^{n}\binom{n}{k}x_ky_{n-k},\quad n=0,1,2,\ldots.$$
\end{DP}
\begin{proof} It is easy to verify that
$z_0z_2=(x_0y_0)(x_0y_2+2x_1y_1+x_2y_0)\ge (x_0y_1+x_1y_0)^2=z_1^2$.
We proceed by induction on $n$. Note that
$$z_n=\sum_{k=0}^{n-1}\binom{n-1}{k}x_ky_{n-k}
+\sum_{k=0}^{n-1}\binom{n-1}{k}x_{k+1}y_{n-k-1}.$$ Two sums in the
right hand side are the binomial convolutions of $\{x_k\}_{0\le k\le
n-1}$ with $\{y_k\}_{1\le k\le n}$ and $\{x_k\}_{1\le k\le n}$ with
$\{y_k\}_{0\le k\le n-1}$ respectively. Hence both are log-convex by
the induction hypothesis. Thus the sequence $\{z_n\}$ is log-convex
by Proposition~\ref{lcx+lcx}.
\end{proof}
\begin{ex}
A permutation $\pi$ of the $n$-element set $[n]=\{1,2,\ldots,n\}$ is
{\it alternating} if $\pi(1)>\pi(2)<\pi(3)>\pi(4)<\cdots \pi(n)$.
The number $E_n$ of alternating permutations of $[n]$ is known as an
{\it Euler number}. The sequence has the exponential generating
function $\sum_{k=0}^{n}E_nx^n/n!=\tan x+\sec x$ and satisfies the
recurrence
$$2E_{n+1}=\sum_{k=0}^{n}\binom{n}{k}E_kE_{n-k}$$
with $E_0=1,E_1=1,E_2=1$ and $E_3=2$ (see
Comtet~\cite[p.~258]{Com74} and Stanley~\cite[p.~149]{Sta97I} for
instance) . Let $z_n=E_n/2$. Then
$z_{n+1}=\sum_{k=0}^{n}\binom{n}{k}z_kz_{n-k}$. Clearly,
$z_0,z_1,z_2,z_3$ is log-convex. Thus the sequence $\{z_n\}$ is
log-convex by induction and Davenport-P\'olya Theorem, and so is the
sequence $\{E_n\}$.
\end{ex}

The following is a special case of Davenport-P\'olya Theorem.
\begin{prop}\label{pascal}
The binomial transformation $z_n=\sum_{k=0}^{n}\binom{n}{k}x_k$
preserves the log-convexity.
\end{prop}

The Stirling number $S(n,k)$ of the second kind is the number of
partitions of the set $[n]$ having exactly $k$ blocks. We have the
following.
\begin{prop}\label{s2-plcx}
The Stirling transformation of the second kind
$z_n=\sum_{k=0}^{n}S(n,k)x_k$ preserves the log-convexity.
\end{prop}
\begin{proof} Let $\{x_k\}_{k\ge 0}$ be a log-convex sequence. We need to
show that the sequence $\{z_n\}_{n\ge 0}$ is log-convex. We proceed
by induction on $n$. It is easy to verify that $z_1^2\le z_0z_2$.
Now assume that $n\ge 3$ and $z_0,z_1,\ldots,z_{n-1}$ is log-convex.
Recall that
$$S(n,k)=\sum_{j=k}^{n}\binom{n-1}{j-1}S(j-1,k-1)$$
(\cite[p.~209]{Com74}). Hence
$$z_n=\sum_{k=1}^{n}\sum_{j=k}^{n}\binom{n-1}{j-1}S(j-1,k-1)x_k
=\sum_{j=0}^{n-1}\binom{n-1}{j}\left[\sum_{k=0}^{j}S(j,k)x_{k+1}\right].$$
Let $y_j=\sum_{k=0}^{j}S(j,k)x_{k+1}$ for $0\le j\le n-1$. Then the
sequence $y_0,y_1,\ldots,y_{n-1}$ is log-convex by the induction
hypothesis, so is the sequence $z_0,z_1,\ldots,z_{n-1},z_{n}$ by
Proposition~\ref{pascal}. This completes the proof. \end{proof}

The Bell number $B_n$ is the total number of partitions of $[n]$,
i.e.,
\begin{equation*}\label{bn}
B_n=\sum_{k=0}^nS(n,k).
\end{equation*}
The log-convexity of Bell numbers was first obtained by
Engel~\cite{Eng94}. Then Bender and Canfield~\cite{BC96} gave a
proof by means of the exponential generating function of the Bell
numbers. Another interesting proof is to use Dobinski formula
(\cite{BT00}). We can also obtain the log-convexity of the Bell
numbers by Proposition~\ref{pascal} and the well-known recurrence
\begin{eqnarray*}
B_{n+1}=\sum_{k=0}^n\binom{n}{k}B_k.
\end{eqnarray*}
Proposition~\ref{s2-plcx} gives a natural interpretation to the
log-convexity of the Bell numbers. The following are some more
applications of Proposition~\ref{s2-plcx}.
\begin{ex}\label{ex-sn}
The Bell number $B_n$ can be viewed as the number of ways of placing
$n$ labeled balls into $n$ indistinguishable boxes. Let $S_n$ be the
number of ways of placing $n$ labeled balls into $n$ unlabeled (but
$2$-colored) boxes. Clearly, $S_n=\sum_{k=1}^{n}2^kS(n,k)$. Hence
the sequence $\{S_n\}_{n\ge 0}$ is log-convex by
Proposition~\ref{s2-plcx}. Note that
$S_n=\sum_{k=0}^{n}\binom{n}{k}B_kB_{n-k}$. Hence the log-convexity
of $\{S_n\}$ can also be followed from Davenport-P\'olya Theorem and
the log-convexity of $\{B_n\}$.
\end{ex}
\begin{ex}\label{ex-cn}
Consider the ordered Bell number $c(n)$, i.e., the number of ordered
partitions of the set $[n]$. Note that $c(n)=\sum_{k=0}^{n}k!S(n,k)$
(Stanley~\cite[p.~146]{Sta97I}). The sequence $\{c(n)\}_{n\ge 0}$ is
therefore log-convex by Proposition~\ref{s2-plcx}.
\end{ex}

Let $c(n,k)$ be the signless Stirling number of the first kind,
i.e., the number of permutations of $[n]$ which contain exactly $k$
permutation cycles. Then
$$c(n,k)=\sum_{j=k}^{n}\binom{n-1}{j-1}(n-j)!c(j-1,k-1)$$
for $1\le k\le n$ (\cite[p.~215]{Com74}). By a similar method used
in the proof of Proposition~\ref{s2-plcx}, we can obtain the
following result.
\begin{prop}
The Stirling transformation of the first kind
$z_n=\sum_{k=0}^{n}c(n,k)x_k$ preserves the log-convexity.
\end{prop}

It is possible to study the log-convexity problems using the theory
of total positivity. Following Karlin~\cite{Kar68}, a matrix
$M=(m_{ij})_{i,j\geq 0}$ of nonnegative numbers is said to be {\it
sign-regular of order r} ($SR_r$ matrix, for short) if, for all
$1\leq t\leq r$, all minors of order $t$ have the same sign.
Denoting with $\varepsilon_{t}\ (=\pm 1)$ the common sign of these
determinants, the vector
$(\varepsilon_{1},\varepsilon_{2},\ldots,\varepsilon_{r})$ is called
the sign sequence of $M$. The matrix is said to be {\it totally
positive of order $r$} (or a $TP_r$ matrix, for short) if it is an
$SR_r$ matrix with $\varepsilon_t=1$ for all $1\leq t\leq r$. A
sequence $a_0,a_1,a_2,\ldots$ of nonnegative numbers is called
$SR_r$ if its Hankel matrix $A=(a_{i+j})_{i,j\geq 0}$ is $SR_r$.
Clearly, the sequence is log-convex (resp. log-concave) if and only
if $a_ma_n\le a_{m-k}a_{n+k}$ (resp. $a_ma_n\ge a_{m-k}a_{n+k}$) for
$1\le k\le m\le n$. In other words, a sequence is log-convex (resp.
log-concave) if and only if it is $SR_2$ with sign sequence $(1,1)$
(resp. $(1,-1)$). As an application of total positivity techniques
to the log-convexity problems, we demonstrate the following
proposition, which is a special case of Brenti~\cite[Theorem
2.2.5]{Bre89}.

Given an infinite lower triangular matrix $A=(a_{n,k})_{n,k\ge 0}$,
let $A_n(u)=\sum_{k=0}^{n}a_{n,k}u^k$ denote the $n$-th row
generating function of $A$.
\begin{prop}\label{lem-bre}
Let $A,B,C$ be three infinite lower triangular matrices satisfying
the following conditions.
\begin{enumerate}
\item [\rm(i)]
Both $B$ and $C$ are $TP_2$ matrices.
\item [\rm(ii)]
$A_{i+j}(u)=B_i(u)C_j(u)$ for all $i,j\in\mathbb{N}$.
\end{enumerate}
Then the linear transformation
\begin{equation*}
z_n=\sum_{k=0}^{n}a_{n,k}x_k,\quad n=0,1,2,\ldots
\end{equation*}
preserves both the log-convexity and the log-concavity.
\end{prop}

In particular, let $P=(p_{n,k})_{n,k\ge 0}$ be the Pascal matrix,
where $p_{n,k}=\binom{n}{k}$ is the binomial coefficients. It is
known that the matrix $P$ is TP$_2$ (all minors of $P$ are actually
nonnegative, see Gessel and Viennot~\cite{GV85} for a combinatorial
proof of this fact). Also, the row generating function
$P_n(u)=(1+u)^n$. Taking $A=B=C=P$ in Proposition~\ref{lem-bre}, we
obtain that the binomial transformation
$z_n=\sum_{k=0}^n\binom{n}{k}x_k$ preserves both the log-convexity
and the log-concavity.
\section{Sequences satisfying three-term recurrences}
\hspace*{\parindent}
In this section we consider the log-convexity of certain famous
combinatorial numbers, including the central binomial coefficients
$b(n)=\binom{2n}{n}$, the Catalan numbers, the Motzkin numbers, the
Fine numbers, the central Delannoy numbers, the little and large
Schr\"oder numbers. These numbers play an important role in
enumerative combinatorics and count various combinatorial objects.
We review some basic facts about these numbers from the viewpoint of
the enumeration of lattice paths in the $(x,y)$ plane.

The central binomial coefficient $b(n)$ counts the number of lattice
paths from $(0,0)$ to $(n,n)$ with steps $(0,1)$ and $(1,0)$ in the
first quadrant. It is clear that $b(n)=\binom{2n}{n}$ for $n\ge 0$.
The central binomial coefficients satisfy the recurrence
$(n+1)b(n+1)=2(2n+1)b(n)$. The sequence $\{b(n)\}_{n\ge 0}$ is
log-convex since $\frac{b(n+1)}{b(n)}=\frac{2(2n+1)}{n+1}$ is
increasing.

The Catalan number $C_n$ counts the number of lattice paths, {\it
Dyck Paths}, from $(0,0)$ to $(2n,0)$ with steps $(1,1)$ and
$(1,-1)$, never falling below the $x$-axis, or equivalently, the
number of lattice paths from $(0,0)$ to $(n,n)$ with steps $(0,1)$
or $(1,0)$, never rising above the line $y=x$. The Catalan numbers
have an explicit expression $C_n=\frac{1}{n+1}\binom{2n}{n}$ and
satisfy the recurrence $(n+2)C_{n+1}=2(2n+1)C_n$. The sequence
$\{C_n\}_{n\ge 0}$ is log-convex since
$\frac{C_{n+1}}{C_n}=\frac{2(2n+1)}{n+2}$ is increasing.

The Motzkin number $M_n$ counts the number of lattice paths, {\it
Motzkin paths}, from $(0,0)$ to $(n,0)$ with steps $(1,0), (1,1)$
and $(1,-1)$, never going below the $x$-axis, or equivalently, the
number of lattice paths from  $(0,0)$ to $(n,n)$, with steps
$(0,2),(2,0)$ and $(1,1)$, never rising above the line $y=x$. It is
known that the Motzkin numbers satisfy the recurrence
\begin{equation}\label{rr-m}
(n+3)M_{n+1}=(2n+3)M_n+3nM_{n-1},
\end{equation}
with $M_0=M_1=1$ (see \cite{Sul01} for a bijective proof).

The Fine number $f_n$ is the number of Dyck paths from $(0,0)$ to
$(2n,0)$ with no hills. (A hill in a Dyck path is a pair of
consecutive steps giving a peak of height $1$). It is known that the
Fine numbers satisfy the recurrence
\begin{equation}\label{rr-f}
2(n+1)f_n=(7n-5)f_{n-1}+2(2n-1)f_{n-2},
\end{equation}
with $f_0=1$ and $f_1=0$ (see \cite{PW99} for a bijective proof).

The central Delannoy number $D(n)$ is the number of lattice paths,
{\it king walks}, from $(0,0)$ to $(n,n)$ with steps $(1,0), (0,1)$
and $(1,1)$ in the first quadrant. Clearly, the number of king walks
with $n-k$ diagonal steps is $\binom{n+k}{n-k}b(k)$. Hence
$D(n)=\sum_{k=0}^{n}\binom{n+k}{n-k}b(k)$ (\cite{Sul03}). It is
known that the central Delannoy numbers satisfy the recurrence
\begin{equation}\label{rr-d}
nD(n)=3(2n-1)D(n-1)-(n-1)D(n-2),
\end{equation}
with $D(0)=1$ and $D(1)=3$ (see \cite{PW02} for a bijective proof).

The (large) Schr\"oder number $r_n$ is the number of king walks,
{\it Schr\"oder paths}, from $(0,0)$ to $(n,n)$, and never rising
above the line $y=x$. The large Schr\"oder numbers bear the same
relation to the Catalan numbers as the central Delannoy numbers do
to the central binomial coefficients. Hence we have
$r_n=\sum_{k=0}^{n}\binom{n+k}{n-k}C_{k}$ (\cite{EHPR98,Wes95}). The
Schr\"oder paths consist of two classes: those with steps on the
main diagonal and those without. These two classes are equinumerous,
and the number of paths in either class is the little Schr\"oder
number $s_n$ (half the large Schr\"oder number). It is known that
the Schr\"oder numbers of two kinds satisfy the recurrence
\begin{equation}\label{rr-ss}
(n+2)z_{n+1}=3(2n+1)z_{n}-(n-1)z_{n-1},
\end{equation}
with $s_0=s_1=r_0=1$ and $r_1=2$ (see Foata and
Zeilberger~\cite{FZ97} for a combinatorial proof and
Sulanke~\cite{Sul98} for another one).

All these numbers presented previously satisfy a three-term
recurrence. Aigner~\cite{Aig98} first established algebraically the
log-convexity of the Motzkin numbers and then Callan~\cite{Cal00}
gave a combinatorial proof. Recently, Do\v sli\'c et
al~\cite{DV03,DSV04} showed the log-convexity of the Motzkin
numbers, the Fine numbers, the central Delannoy numbers, the large
and little Schr\"oder numbers using calculus. (We appreciate Do\v
sli\'c for acquainting us with his very recent work
\cite{Dos05,Dos06}, in which the techniques developed in
\cite{DV03,DSV04} are extended and more examples are presented.)
Motivated by these results, we investigate the log-convexity problem
of combinatorial sequences satisfying a three-term recurrence by an
algebraic approach. We distinguish two cases according to the sign
of coefficients in the recurrence relations.
\subsection{The recurrence $a_nz_{n+1}=b_nz_{n}+c_nz_{n-1}$}
\hspace*{\parindent}
Let $\{z_n\}_{n\ge 0}$ be a sequence of positive numbers satisfying
the recurrence
\begin{equation}\label{rr-c+}
a_nz_{n+1}=b_nz_{n}+c_nz_{n-1}
\end{equation}
for $n\ge 1$, where $a_n,b_n,c_n$ are all positive. Consider the
quadratic equation
\begin{equation*}
a_n\la^2-b_n\la-c_n=0
\end{equation*}
associated with the recurrence~(\ref{rr-c+}).
Clearly, the equation has a unique positive root
\begin{equation}\label{lambda}
\la_n=\frac{b_n+\sqrt{b_n^2+4a_nc_n}}{2a_n}.
\end{equation}
Define $x_n=z_{n+1}/z_{n}$ for $n\ge 0$. Then the sequence
$\{z_n\}_{n\ge 0}$ is log-convex if and only if the sequence
$\{x_n\}_{n\ge 0}$ is increasing. By (\ref{rr-c+}), the sequence
$\{x_n\}_{n\ge 0}$ satisfies the recurrence
\begin{equation}\label{rr-x}
a_nx_n=b_n+\frac{c_n}{x_{n-1}}
\end{equation}
for $n\ge 1$. It follows that $x_{n}\ge x_{n-1}$ is equivalent to
$x_{n-1}\le\la_n$, and is also equivalent to $x_n\ge\la_n$. Thus the
sequence $\{z_n\}_{n\ge 0}$ is log-convex if and only if the
sequence $\{x_n\}_{n\ge 0}$ can be separated by the sequence
$\{\la_n\}_{n\ge 1}$:
\begin{equation}\label{int}
x_0\le\la_1\le x_1\le\la_2\le\cdots\le x_{n-1}\le\la_n\le
x_n\le\la_{n+1}\le\cdots.
\end{equation}
\begin{thm}\label{crit+}
Let $\{z_n\}_{n\ge 0}$ and $\{\la_n\}_{n\ge 1}$ be as above. Suppose
that $z_0,z_1,z_2,z_3$ is log-convex and that the inequality
\begin{equation}\label{la-cond}
a_n\la_{n-1}\la_{n+1}-b_n\la_{n-1}-c_n\ge 0 \end{equation} holds for
$n\ge 2$. Then the sequence $\{z_n\}_{n\ge 0}$ is log-convex.
\end{thm}
\begin{proof} Let $x_n=z_{n+1}/z_n$ for $n\ge 0$. We prove the interlacing
inequalities (\ref{int}) by induction. By the assumption that
$z_0,z_1,z_2,z_3$ is log-convex, we have 
$x_0\le\la_1\le x_1\le\la_2\le x_2$. Now assume that $\la_{n-1}\le
x_{n-1}\le\la_{n}$. Note that $x_{n-1}\le\la_n$ is equivalent to
$x_n\ge\la_n$. On the other hand, $x_{n-1}\ge\la_{n-1}$ implies that
$$x_n=\frac{b_n}{a_n}+\frac{c_n}{a_n}\frac{1}{x_{n-1}}\le \frac{b_n}{a_n}+\frac{c_n}{a_n}\frac{1}{\la_{n-1}}\le\la_{n+1}$$
by the inequality~(\ref{la-cond}). Hence we have $\la_{n}\le
x_{n}\le\la_{n+1}$. Thus (\ref{int}) holds by induction. \end{proof}
\begin{cor}
The Fine sequence $\{f_n\}_{n\ge 2}$ is log-convex.
\end{cor}
\begin{proof} Let $z_n=f_{n+2}$ for $n\ge 0$. Then
$z_0=1,z_1=2,z_2=6,z_3=18$ and
$$2(n+4)z_{n+1}=(7n+16)z_n+2(2n+5)z_{n-1}$$
by (\ref{rr-f}). 
Solve the equation $2(n+4)\la^2-(7n+16)\la-2(2n+5)=0$ to obtain
$$\la_n=\frac{2(2n+5)}{(n+4)}.$$ It is easy to verify that
$$2(n+4)\la_{n-1}\la_{n+1}-(7n+6)\la_{n-1}-2(2n+5)=\frac{2(20n^2+151n+171)}{(n+3)(n+5)}\ge 0.$$
Hence the sequence $\{z_n\}_{n\ge 0}$, i.e., $\{f_n\}_{n\ge 2}$, is
log-convex by Theorem~\ref{crit+}.\end{proof}

Because of the expression (\ref{lambda}) of $\la_n$, sometimes it is
inconvenient to directly check the inequality~(\ref{la-cond}).
However, the inequality can be verified by means of Maple. For
example, for the Motzkin sequence satisfying the
recurrence~(\ref{rr-m}), we have
$$\la_n=\frac{2n+3+\sqrt{16n^2+48n+9}}{2(n+3)}.$$ Using Maple it is
easy to verify the inequality
$$(n+3)\la_{n-1}\la_{n+1}-(2n+3)\la_{n-1}-3n\ge 0.$$
Thus the log-convexity of the Motzkin numbers follows from
Theorem~\ref{crit+}.
\begin{cor}
The Motzkin sequence $\{M_n\}_{n\ge 0}$ is log-convex.
\end{cor}

We can also give another criterion for the log-convexity of the
sequence $\{z_n\}$ satisfying the recurrence (\ref{rr-c+}).
\begin{thm}\label{c+}
Let $\{z_n\}_{n\ge 0}$ and $\{\la_n\}_{n\ge 1}$ be defined by
(\ref{rr-c+}) and (\ref{lambda}). Suppose that there exists a
sequence $\{\mu_n\}_{n\ge 1}$ of positive numbers such that the
following three conditions hold.
\begin{enumerate}
\item [\rm(i)]
$\mu_n\le\la_n$ for all $n\ge 1$.
\item [\rm(ii)]
$z_1\le \mu_1z_0$ and $z_2\le\mu_2z_1$.
\item [\rm(iii)]
$a_n\mu_{n-1}\mu_{n+1}\ge b_n\mu_{n-1}+c_n$ for $n\ge 2$.
\end{enumerate}
Then the sequence $\{z_n\}_{n\ge 0}$ is log-convex.
\end{thm}
\begin{proof}
Let $x_n=z_{n+1}/z_n$ for $n\ge 0$. Then it suffices to show
that the sequence $\{x_n\}$ is increasing. We prove this by showing
the interlacing inequalities
\begin{equation}\label{int-mu}
x_0\le\mu_1\le x_1\le\mu_2\le\cdots\le x_{n-1}\le\mu_n\le
x_n\le\mu_{n+1}\le\cdots.
\end{equation}
The condition (ii) is equivalent to $x_0\le\mu_1$ and $x_1\le\mu_2$.
However, $\mu_1\le\la_1$. Hence $x_0\le\la_1$, and so
$x_1\ge\la_1\ge\mu_1$. Thus we have $\mu_1\le x_1\le\mu_2$. Now
assume that $\mu_{n-1}\le x_{n-1}\le\mu_n$. Then $x_{n-1}\le\mu_n$
implies $x_n\ge\mu_n$ since $\mu_n\le\la_n$. On the other hand,
$x_{n-1}\ge\mu_{n-1}$ implies
$$x_n=\frac{b_n}{a_n}+\frac{c_n}{a_n}\frac{1}{x_{n-1}}
\le\frac{b_n}{a_n}+\frac{c_n}{a_n}\frac{1}{\mu_{n-1}}\le\mu_{n+1}$$
by the condition (iii). Hence we have $\mu_n\le x_n\le\mu_{n+1}$.
Thus (\ref{int-mu}) holds by induction.
\end{proof}

For convenience, we may choose $\mu_n$ in the theorem as an
appropriate rational approximation to $\la_n$. We present two
examples to demonstrate this approach.

The derangements number $d_n$ is the number of permutations of $n$
elements with no fixed points. It is well known that the sequence
$\{d_n\}_{n\ge 0}$ satisfies the recurrence
$$d_{n+1}=n(d_{n}+d_{n-1}),$$ with $d_0=1,d_1=0,d_2=1,d_3=2$ and
$d_4=9$ (Comtet~\cite[p.~182]{Com74}).
\begin{cor}
The sequence of the derangements numbers $\{d_n\}_{n\ge 2}$ is
log-convex.
\end{cor}
\begin{proof}
Let $z_n=d_{n+2}$ for $n\ge 0$. Then the sequence satisfies
the recurrence
$$z_{n+1}=(n+2)(z_{n}+z_{n-1}),$$
with $z_0=1,z_1=2$ and $z_2=9$. We have
$$\la_n=\frac{(n+2)+\sqrt{n^2+8n+12}}{2}\ge \frac{(n+2)+(n+3)}{2}=\frac{2n+5}{2}.$$
Set $\mu_n=(2n+5)/2$. Then $\mu_n\le\la_n$. Also,
$z_1/z_0=2<7/2=\mu_1$ and $z_2/z_1=9/2=\mu_2$. Furthermore,
$$\mu_{n-1}\mu_{n+1}-(n+2)\mu_{n-1}-(n+2)
=\frac{2n+1}{4}\ge 0.$$ Thus the sequence $\{z_n\}_{n\ge 0}$, i.e.,
$\{d_n\}_{n\ge 2}$, is log-convex by Theorem~\ref{c+}.
\end{proof}
\begin{rem}
Generally, if $a_n$ takes a constant value, both $b_n$ and $c_n$ are
linear functions in $n$ respectively, then we can show that the
sequence $\{z_n\}_{n\ge 0}$ satisfying the recurrence~(\ref{rr-c+})
is asymptotically log-convex by means of Theorem~\ref{crit+}. In
other words, there exists an index $N$ such that $\{z_n\}_{n\ge N}$
is log-convex. We leave the proof of this result to the reader as an
exercise.
\end{rem}

Let $A_n$ be the number of directed animals of size $n$
(Stanley~\cite[Exercise 6.46]{Sta97II}). The sequence $\{A_n\}_{n\ge
0}$ is Sloane's A005773 (\cite{Slo}) and satisfies the recurrence
$$(n+1)A_{n+1}=2(n+1)A_{n}+3(n-1)A_{n-1},$$ with $A_0=1,A_1=1$ and
$A_2=2$.
\begin{cor}
The sequence $\{A_n\}_{n\ge 0}$ counting directed animals is
log-convex.
\end{cor}
\begin{proof}
Note that
$$\la_n=1+\sqrt{\frac{4n-2}{n+1}}=1+\sqrt{\frac{4n-2}{2n}\frac{4n}{2n+2}}\ge
1+\frac{4n-1}{2n+1}=\frac{6n}{2n+1}$$ for $n\ge 1$, where the
inequality follows from
$\sqrt{\frac{x-1}{y-1}\frac{x+1}{y+1}}\ge\frac{x}{y}$ when $x\ge
y>1$. Set $\mu_n=\frac{6n}{2n+1}$. Then $\mu_n\le\la_n,
A_1/A_0=1<2=\mu_1$ and $A_2/A_1=2<2.4=\mu_2$. Also,
$$(n+1)\mu_{n-1}\mu_{n+1}-2(n+1)\mu_{n-1}-3(n-1)=\frac{9(n-1)}{(2n-1)(2n+3)}\ge 0.$$
Thus the sequence $\{A_n\}_{n\geq 0}$ is log-convex by
Theorem~\ref{c+}.
\end{proof}

The techniques developed in this subsection can also be used to
study the log-concavity of the sequences $\{z_n\}_{n\ge 0}$
satisfying the recurrence~(\ref{rr-c+}). For example, it is clear
that $\{z_n\}_{n\ge 0}$ is log-concave if and only if
\begin{equation*}\label{int-lc}
x_0\ge\la_1\ge x_1\ge\la_2\ge\cdots\ge
x_{n-1}\ge\la_n\ge x_n\ge\la_{n+1}\ge\cdots.
\end{equation*}
Then we have the following result similar to Theorem~\ref{crit+}.
\begin{thm}\label{lc+} Let $\{z_n\}_{n\ge 0}$ and $\{\la_n\}_{n\ge 1}$ be defined by (\ref{rr-c+}) and (\ref{lambda}).
Suppose that $z_0,z_1,z_2,z_3$ is log-concave and that the
inequality
\begin{equation*}
a_n\la_{n-1}\la_{n+1}-b_n\la_{n-1}-c_n\le 0 \end{equation*} holds
for $n\ge 2$. Then the sequence $\{z_n\}_{n\ge 0}$ is log-concave.
\end{thm}
\begin{rem}
When coefficients $a_n,b_n,c_n$ in the recurrence~(\ref{rr-c+}) take
constant values respectively, the sequence $\{z_n\}_{n\ge 0}$ is
neither log-convex nor log-concave since $\la_n$ takes a constant
value. For example, the Fibonacci numbers satisfy the recurrence
$F_{n+1}=F_n+F_{n-1}$ with $F_0=F_1=1$. The sequence $\{F_n\}_{n\ge
0}$ is neither log-convex nor log-concave. Actually,
$F_{n-1}F_{n+1}-F^2_n=(-1)^{n-1}$. However, the bisection
$\{F_{2n}\}_{n\ge 0}$ with even index is log-convex and the
bisection $\{F_{2n+1}\}_{n\ge 0}$ with odd index is log-concave
since $F_{n-2}F_{n+2}-F^2_n=(-1)^{n}$. We will give a general result
about such sequences in Corollary~\ref{bisection}.
\end{rem}
\subsection{The recurrence $a_nz_{n+1}=b_nz_{n}-c_nz_{n-1}$}
\hspace*{\parindent}
In this part we consider the log-convexity of the sequence $\{z_n\}$
of positive numbers satisfying the recurrence
\begin{equation}\label{rr-c-}
a_nz_{n+1}=b_nz_{n}-c_nz_{n-1}
\end{equation}
for $n\ge 1$, where $a_n,b_n,c_n$ are all positive.

 Let $x_n=z_{n+1}/z_{n}$ for $n\ge 0$. Then we need to check whether the
sequence $\{x_n\}_{n\ge 0}$ is increasing.

By the recurrence~(\ref{rr-c-}), we have
\begin{equation*}\label{x:n-1->n}
x_n=\frac{b_n}{a_n}-\frac{c_n}{a_n}\frac{1}{x_{n-1}}.
\end{equation*}
Hence
\begin{equation}\label{x:n+1-n}
x_{n+1}-x_n=\left[\left(\frac{b_{n+1}}{a_{n+1}}-\frac{b_n}{a_n}\right)
+\left(\frac{c_n}{a_n}-\frac{c_{n+1}}{a_{n+1}}\right)\frac{1}{x_{n}}\right]
+\frac{c_n}{a_n}\left(\frac{1}{x_{n-1}}-\frac{1}{x_n}\right).
\end{equation}
Observe that if
\begin{equation}\label{nn+1}
\begin{vmatrix}a_n & a_{n+1}\\
b_n & b_{n+1}\end{vmatrix}x_n +
\begin{vmatrix}c_n & c_{n+1}\\
a_n & a_{n+1}\end{vmatrix} \ge 0,
\end{equation}
then $x_{n-1}\le x_n$ implies $x_n\le x_{n+1}$ from (\ref{x:n+1-n}).
Hence we can conclude that if $x_0\le x_1$ and the inequality
(\ref{nn+1}) holds for $n\ge 1$, then the sequence $\{x_n\}_{n\ge
0}$ is increasing, and the sequence $\{z_n\}$ is therefore
log-convex.

Suppose now that $a_n,b_n,c_n$ are all linear functions in $n$. In
this case, the inequality (\ref{nn+1}) is easily checked since two
determinants take constant values respectively. More precisely, let
$$a_n=\alpha_1n+\alpha_0,\quad b_n=\beta_1n+\beta_0,\quad
c_n=\gamma_1n+\gamma_0$$ and denote
\begin{equation*}\label{abc}A=\begin{vmatrix}\beta_0 & \beta_1\\ \gamma_0 &
\gamma_1\end{vmatrix},\quad B=\begin{vmatrix}\gamma_0 &
\gamma_1\\\alpha_0 & \alpha_1
\end{vmatrix},\quad C=\begin{vmatrix}\alpha_0 & \alpha_1\\
\beta_0 & \beta_1\end{vmatrix}.\end{equation*} Then it is easy to
see that two determinants in the inequality (\ref{nn+1}) are equal
to $C$ and $B$ respectively. Thus we have the following criterion.
\begin{thm}\label{c-}
Let $\{z_n\}_{n\ge 0}$ be a sequence of positive numbers and satisfy
the three-term recurrence
\begin{equation}\label{lin-}
(\alpha_1n+\alpha_0)z_{n+1}=(\beta_1n+\beta_0)z_{n}-(\gamma_1n+\gamma_0)z_{n-1}
\end{equation}
for $n\ge 1$, where
$\alpha_1n+\alpha_0,\beta_1n+\beta_0,\gamma_1n+\gamma_0$ are
positive for $n\ge 1$. Suppose that $z_0,z_1,z_2$ is log-convex.
Then the full sequence $\{z_n\}_{n\ge 0}$ is log-convex if one of
the following conditions holds.
\begin{itemize}
  \item [\rm(i)] $B,C\ge 0$.
  \item [\rm(ii)] $B<0,C>0,AC\ge B^2$ and $z_0B+z_1C\ge 0$.
  \item [\rm(iii)] $B>0,C<0,AC\le B^2$ and $z_0B+z_1C\ge 0$.
\end{itemize}
\end{thm}
\begin{proof}
Let $x_n=z_{n+1}/z_n$ for $n\ge 0$. Then $x_0\le x_1$ since
$z_0z_2\ge z_1^2$. Thus it suffices to show that the inequality
$Cx_n+B\ge 0$ holds for $n\ge 0$. If $B\ge 0$ and $C\ge 0$, then the
inequality is obvious. Next we assume that $BC<0$ and show that
$Cx_n+B\ge 0$ by induction on $n$. We do it only for the case (ii)
since the case (iii) is similar. Clearly, $Cx_0+B\ge 0$ by the
condition $z_0B+z_1C\ge 0$. Now assume that $Cx_{n-1}+B\ge 0$ for
$n\ge 1$. Then
\begin{eqnarray*}
Cx_n+B &=& C\left(\frac{b_n}{a_n}-\frac{c_n}{a_n}\frac{1}{x_{n-1}}\right)+B\\
&\ge & \frac{C}{a_n}\left(b_n+c_n\frac{C}{B}\right)+B\\
& =& \frac{C}{a_nB}(b_nB+c_nC)+B.
\end{eqnarray*}
Note that $b_nB+c_nC=-a_nA$ since
$$a_nA+b_nB+c_nC=\begin{vmatrix}
a_n & \alpha_1 & \alpha_0\\
b_n & \beta_1 & \beta_0\\
c_n & \gamma_1 & \gamma_0 \end{vmatrix}=
\begin{vmatrix}
\alpha_1n+\alpha_0 & \alpha_1 & \alpha_0\\
\beta_1n+\beta_0 & \beta_1 & \beta_0\\
\gamma_1n+\gamma_0 & \gamma_1 & \gamma_0 \end{vmatrix} =0.$$ Hence
$$Cx_n+B\ge -\frac{AC}{B}+B\ge 0$$ by the condition $AC\ge B^2$, as desired.
This completes our proof.
\end{proof}

The techniques developed in Theorem~\ref{c-} can also be used to
study the log-concavity of the sequences $\{z_n\}_{n\ge 0}$
satisfying the recurrence~(\ref{lin-}). We demonstrate the result
without proof.
\begin{thm}\label{rr-c-lc}
Let $\{z_n\}_{n\ge 0}$ be a sequence of positive numbers and satisfy
the three-term recurrence (\ref{lin-}). Suppose that $z_0,z_1,z_2$
is log-concave. Then the full sequence $\{z_n\}_{n\ge 0}$ is
log-concave if one of the following conditions holds.
\begin{itemize}
 \item [\rm(i)] $B,C\le 0$.
  \item [\rm(ii)] $B<0,C>0,AC\le B^2$ and $z_0B+z_1C\le 0$.
  \item [\rm(iii)] $B>0,C<0,AC\ge B^2$ and $z_0B+z_1C\le 0$.
\end{itemize}
\end{thm}
\begin{cor}
The central Delannoy sequence $\{D(n)\}_{n\ge 0}$ is log-convex.
\end{cor}
\begin{proof}
By the recurrence (\ref{rr-d}), we have $A=3, B=-1, C=3$. Also,
$D(0)=1, D(1)=3, D(2)=13$. Thus the log-convexity of $\{D(n)\}_{n\ge
0}$ follows from Theorem~\ref{c-}.
\end{proof}
\begin{cor}
The little and large Schr\"oder numbers are log-convex respectively.
\end{cor}
\begin{proof}
It suffices to show that the little Schr\"oder numbers
$\{s_n\}_{n\ge 0}$ is log-convex since the large Schr\"oder numbers
$r_n=2s_n$ for $n\ge 1$.

By the recurrence (\ref{rr-ss}), we have $A=9,B=-3,C=9$. Also,
$s_0=s_1=1,s_2=3$. Thus the log-convexity of $\{s_n\}_{n\ge 0}$
follows from Theorem~\ref{c-}.
\end{proof}

Let $h_n$ be the number of the set of all tree-like polyhexes with
$n+1$ hexagons (Harary and Read~\cite{HR70}). It is known that $h_n$
counts the number of lattice paths, from $(0,0)$ to $(2n,0)$ with
steps $(1,1),(1,-1)$ and $(2,0)$, never falling below the $x$-axis
and with no peaks at odd level. The sequence $\{h_n\}_{n\ge 0}$ is
Sloane's A002212 and satisfies the recurrence
$$(n+1)h_n=3(2n-1)h_{n-1}-5(n-2)h_{n-2}$$with $h_0=h_1=1$ and $h_2=3$.
Thus the following corollary is an immediate consequence of Theorem
3.10.
\begin{cor}
The sequence $\{h_n\}_{n\ge 0}$ counting tree-like polyhexes is
log-convex.
\end{cor}

Let $w_n$ be the number of walks on cubic lattice with $n$ steps,
starting and finishing on the $xy$ plane and never going below it
(Guy~\cite{Guy00}). The sequence $\{w_n\}_{n\ge 0}$ is Sloane's
A005572 and satisfies the recurrence
$$(n+2)w_n=4(2n+1)w_{n-1}-12(n-1)w_{n-2},$$
with $w_0=1, w_1=4$ and $w_2=17$. Thus the following corollary is
immediate from Theorem 3.10.
\begin{cor}
The sequence $\{w_n\}_{n\ge 0}$ counting walks on cubic lattice is
log-convex.
\end{cor}

A special interesting case of Theorem~\ref{c-} (i) and
Theorem~\ref{rr-c-lc} (i) is the following.
\begin{cor}\label{c-c}
Suppose that the sequence $\{z_n\}_{n\ge 0}$ of positive numbers
satisfies the recurrence $az_{n+1}=bz_n-cz_{n-1}$ for $n\ge 1$,
where $a,b,c$ are positive constants. If $z_0,z_1,z_2$ is log-convex
(resp. log-concave), then so is the full sequence $\{z_n\}_{n\ge
0}$.
\end{cor}
\begin{cor}\label{bisection}
Suppose that the sequence $\{z_n\}_{n\ge 0}$ of positive numbers
satisfies the recurrence $az_{n+1}=bz_n+cz_{n-1}$ for $n\ge 1$,
where $a,b,c$ are positive constants. If $z_0,z_1,z_2$ is log-convex
(resp. log-concave), then the bisection $\{z_{2n}\}$ is log-convex
(resp. log-concave) and the bisection $\{z_{2n+1}\}$ is log-concave
(resp. log-convex).
\end{cor}
\begin{proof}
By the recurrence $az_{n+1}=bz_n+cz_{n-1}$ for $n\ge 1$, we
can obtain the recurrence
$$a^2z_{n+2}=(b^2+2ac)z_n-c^2z_{n-2}$$
for $n\ge 2$. It is not difficult to verify that
$$a^2(z_0z_4-z^2_2)=b^2(z_0z_2-z^2_1),
\quad a^3(z_1z_5-z^2_3)=b^2c(z_1^2-z_0z_2).$$ So the statement
follows from Corollary~\ref{c-c}.
\end{proof}

The Fibonacci numbers $F_n$ satisfy the recurrence
$F_{n+1}=F_n+F_{n-1}$ with $F_0=F_1=1$ and $F_2=2$. The Lucas
numbers $L_n$ satisfy the recurrence $L_{n+1}=L_n+L_{n-1}$ with
$L_0=1,L_1=3$ and $L_2=4$. And the Pell numbers $P_n$ satisfy the
recurrence $P_{n+1}=2P_n+P_{n-1}$ with $P_0=1, P_1=2$ and $P_3=5$.
Thus we can conclude the following result from
Corollary~\ref{bisection}.
\begin{cor}\label{rr-c-+}
The bisections $\{F_{2n+1}\}$, $\{L_{2n}\}$, $\{P_{2n+1}\}$ are
log-concave and the bisections $\{F_{2n}\}$, $\{L_{2n+1}\}$,
$\{P_{2n}\}$ are log-convex.
\end{cor}
\section{$q$-log-convexity}
\hspace*{\parindent}
In this section we first introduce the concept of the
$q$-log-convexity of polynomial sequences and then prove the
$q$-log-convexity of certain well-known polynomial sequences,
including the Bell polynomials, the Eulerian polynomials, the
$q$-Schr\"oder numbers and the $q$-central Delannoy numbers. We also
present certain linear transformations preserving the log-convexity
of sequences and establish the connection with the
$q$-log-convexity.

Let $q$ be an indeterminate. Given two real polynomials $f(q)$ and
$g(q)$, write $f(q)\le_q g(q)$ if and only if $g(q)-f(q)$ has
nonnegative coefficients as a polynomial in $q$. A sequence of real
polynomials $\{P_n(q)\}_{n\ge 0}$ is called {\it $q$-log-convex} if
\begin{equation}\label{q-lcx}
P^2_n(q)\le_q P_{n-1}(q)P_{n+1}(q)
\end{equation}
for all $n\ge 1$. Clearly, if the sequence $\{P_n(q)\}_{n\ge 0}$ is
$q$-log-convex, then for each fixed positive number $q$, the
sequence $\{P_n(q)\}_{n\ge 0}$ is log-convex. The converse is not
true in general. If the opposite inequality in (\ref{q-lcx}) holds,
then the sequence $\{P_n(q)\}_{n\ge 0}$ is called {\it
$q$-log-concave}. The concept of the $q$-log-concavity was first
suggested by Stanley and these has been much interest in this
subject. We refer the reader to Sagan~\cite{Sag92TAMS,Sag92DM} for
further information about the $q$-log-concavity.

Perhaps the simplest example of $q$-log-convex polynomials is the
$q$-factorial. It is well known that the factorial $n!$ is
log-convex. The standard $q$-analogue of an integer $n$ is
$(n)_q=1+q+q^2+\cdots+q^{n-1}$ and the associated $q$-factorial is
${(n)_q}!=\prod_{k=1}^{n}(k)_q$. It is easy to verify that the
$q$-factorial ${(n)_q}!$ is $q$-log-convex by a direct calculation.
We next provide more examples of $q$-log-convex sequences.
\subsection{Bell polynomials and Eulerian polynomials}
\hspace*{\parindent}
The Bell polynomial, or the exponential polynomial, is the
generating function $B_n(q)=\sum_{k=0}^{n}S(n,k)q^k$ of the Stirling
numbers of the second kind. It can be viewed as a $q$-analog of the
Bell number and has many fascinating properties (see Roman~\cite[\S
4.1.3]{Rom84} for instance). Note that the Stirling numbers of the
second kind satisfy the recurrence
\begin{equation*}
S(n+1,k)=kS(n,k)+S(n,k-1)
\end{equation*}
Hence the Bell polynomials satisfy the recurrence
$$B_{n+1}(q)=qB_n(q)+qB'_n(q).$$
It is well known that the Bell polynomials $B_n(q)$ have only real
zeros (see \cite{WYjcta05} for instance). In \S2 we have shown that
the linear transformation $z_n=\sum_{k=0}^{n}S(n,k)x_k$ can preserve
the log-convexity of sequences. Therefore, for each positive number
$q$, the sequence $\{B_n(q)\}_{n\ge 0}$ is log-convex. A further
problem is whether the sequence $\{B_n(q)\}_{n\ge 0}$ is
$q$-log-convex.

Let $\pi=a_1a_2\cdots a_n$ be a permutation of $[n]$. An element
$i\in [n-1]$ is called a descent of $\pi$ if $a_i>a_{i+1}$. The
Eulerian number $A(n,k)$ is defined as the number of permutations of
$[n]$ having $k-1$ descents.
 The Eulerian numbers satisfy the recurrence
\begin{equation*}\label{rr-an}
A(n,k)=kA(n-1,k)+(n-k+1)A(n-1,k-1).
\end{equation*}
Let $A_n(q)=\sum_{k=0}^{n}A(n,k)q^k$ be the Eulerian polynomial.
Then $$A_n(q)=nqA_{n-1}(q)+q(1-q)A'_{n-1}(q).$$ It is well known
that $A_n(q)$ has only real zeros and $A(n,k)$ is therefore
log-concave in $k$ for each fixed $n$ (see \cite{WYjcta05} for
instance). By Frobenius formula
$$A_n(q)=q\sum_{k=1}^{n}k!S(n,k)(q-1)^{n-k}$$ and Proposition~\ref{s2-plcx}, the sequence $\{A_n(q)\}_{n\ge 0}$ is
log-convex for each fixed positive number $q\ge 1$. We refer the
reader to Comtet~\cite{Com74} for further information about the
Eulerian numbers and the Eulerian polynomials.

To show the $q$-log-convexity of both the Bell polynomials and the
Eulerian polynomials, we establish the following more general
result.

Let $\{T(n,k)\}_{n,k\ge 0}$ be an array of nonnegative numbers
satisfying the recurrence
\begin{equation}\label{rr-tnk}
T(n,k)=(a_1n+a_2k+a_3)T(n-1,k)+(b_1n+b_2k+b_3)T(n-1,k-1)
\end{equation}
with $T(n,k)=0$ unless $0\le k\le n$. It is natural to assume that
$a_1n+a_2k+a_3\ge 0$ for $0\le k<n$ and $b_1n+b_2k+b_3\ge 0$ for
$0<k\le n$. Note that the former is equivalent to $a_1\ge
0,a_1+a_2\ge 0, a_1+a_3\ge 0$ and the latter is equivalent to
$b_1\ge 0,b_1+b_2\ge 0,b_1+b_2+b_3\ge 0$. It is known that for each
fixed $n$, the sequence $\{T(n,k)\}_{0\le k\le n}$ is log-concave
(Kurtz~\cite{Kur72}) and further, is PF if $a_2b_1\ge a_1b_2$ and
$a_2(b_1+b_2+b_3)\ge (a_1+a_3)b_2$ (Wang and Yeh~\cite[Corollary
3]{WYjcta05}).
\begin{thm}\label{T-qLCX}
Let $\{T(n,k)\}_{n,k\ge 0}$ be as above and the row generating
function $T_n(q)=\sum_{k=0}^{n}T(n,k)q^k$. Suppose that for $0<k\le
n$,
\begin{eqnarray}\label{T-condition}
(a_2b_1-a_1b_2)n+a_2b_2k+(a_2b_3-a_3b_2)\ge 0.
\end{eqnarray}
Then the sequence $\{T_n(q)\}_{n\ge 0}$ is $q$-log-convex.
\end{thm}
\begin{rem}
The condition (\ref{T-condition}) is equivalent to
$$a_2b_1-a_1b_2, a_2(b_1+b_2)-a_1b_2,
a_2(b_1+b_2+b_3)-(a_1+a_3)b_2$$ are all nonnegative. Hence the
polynomial $T_n(q)$ in Theorem~\ref{T-qLCX} has only real zeros for
each $n\ge 0$ by Wang and Yeh~\cite[Corollary 3]{WYjcta05}.
\end{rem}
\begin{rem}
If $a_2=b_2=0$, then the condition (\ref{T-condition}) is trivially
satisfied.
\end{rem}
\begin{proof}[Proof of Theorem~\ref{T-qLCX}] Let
$T_{n-1}(q)T_{n+1}(q)-T_n^2(q)=\sum_{t=0}^{2n}A_tq^t$. We need to
show that $A_t\ge 0$ for $0\le t\le 2n$. Note that the recurrence
(\ref{rr-tnk}) is equivalent to
\begin{eqnarray*}
T_n(q) &=&
(a_1n+a_3+b_1nq+b_2q+b_3q)T_{n-1}(q)+(a_2+b_2q)qT'_{n-1}(q).
\end{eqnarray*}
Hence
\begin{eqnarray*}
\sum_{t=0}^{2n}A_tq^t
&=& T_{n-1}(q)[(a_1n+a_1+a_3+b_1nq+b_1q+b_2q+b_3q)T_{n}(q)+(a_2+b_2q)qT'_{n}(q)]\\
& & -T_n(q)[(a_1n+a_3+b_1nq+b_2q+b_3q)T_{n-1}(q)+(a_2+b_2q)qT'_{n-1}(q)]\\
&=&
(a_1+b_1q)T_{n-1}(q)T_n(q)+(a_2+b_2q)q[T_{n-1}(q)T'_n(q)-T'_{n-1}(q)T_n(q)].
\end{eqnarray*}
Thus $A_t=\sum_{k=0}^{t}c_k(n,t)$, where
\begin{eqnarray*}
c_k
&=& T(n,t-k)[a_1T(n-1,k)+b_1T(n-1,k-1)+a_2(t-k)T(n-1,k)\\
& & -a_2kT(n-1,k)+b_2(t-k)T(n-1,k-1)-b_2(k-1)T(n-1,k-1)]\\
&=&
T(n,t-k)[(a_1+ta_2-2ka_2)T(n-1,k)+(b_1+b_2+tb_2-2kb_2)T(n-1,k-1)].
\end{eqnarray*}
Clearly, $c_k\ge 0$ if $t$ is even and $k=t/2$. So, in order to
prove that $A_t\ge 0$, it suffices to prove that $c_k+c_{t-k}\ge 0$
for $0\le k<t-k\le n$. Let $u_k=T(n-1,k)$ if $0\le k\le n-1$ and
$u_k=0$ otherwise. Then the sequence $\{u_k\}_{k\ge 0}$ is
log-concave. In what follows, we always assume that $0\le k<t-k\le
n$. By the recurrence (\ref{rr-tnk}), we have
\begin{eqnarray*}
c_k+c_{t-k}
&=& [(a_1n+a_2(t-k)+a_3)u_{t-k}+(b_1n+b_2(t-k)+b_3)u_{t-k-1}]\\
& & \times[(a_1+ta_2-2ka_2)u_k+(b_1+b_2+tb_2-2kb_2)u_{k-1}]\\
& & +[(a_1n+a_2k+a_3)u_{k}+(b_1n+b_2k+b_3)u_{k-1}]\\
& & \times[(a_1-ta_2+2ka_2)u_{t-k}+(b_1+b_2-tb_2+2kb_2)u_{t-k-1}]\\
&=& xu_ku_{t-k}+yu_{k-1}u_{t-k-1}+zu_{k}u_{t-k-1}+wu_{k-1}u_{t-k},
\end{eqnarray*}
where
\begin{eqnarray*}
x&=& a_1(2a_1n+a_2t+2a_3)+a_2^2(t-2k)^2,\\
y&=& (b_1+b_2)(2b_1n+b_2t+2b_3)+b_2^2(t-2k)^2,\\
z&=& (a_1n+a_2k+a_3)(b_1+b_2)+[b_1n+b_2(t-k)+b_3]a_1\\
& & +(t-2k)[(a_2b_1-a_1b_2)n+a_2b_2(t-2k)+(a_2b_3-a_3b_2)],\\
w&=& (b_1n+b_2k+b_3)a_1+[a_1n+a_2(t-k)+a_3](b_1+b_2)\\
& & -(t-2k)[(a_2b_1-a_1b_2)n-a_2b_2(t-2k)+(a_2b_3-a_3b_2)].
\end{eqnarray*}
By the assumption of the recurrence (\ref{rr-tnk}), the numbers
$a_1,b_1+b_2$ and $a_1n+a_2k+a_3,b_1n+b_2(t-k)+b_3$ are all
nonnegative for $0\le k<t-k\le n$. On the other hand, by the
assumption of the theorem, the number
$(a_2b_1-a_1b_2)n+a_2b_2(t-2k)+(a_2b_3-a_3b_2)$ is nonnegative since
$0<t-2k\le t-k\le n$. Hence $z\ge 0$ for $0\le k<t-k\le n$. Thus by
the log-concavity of $\{u_k\}$, we have
\begin{eqnarray*}
c_k+c_{t-k}
&\ge & xu_ku_{t-k}+yu_{k-1}u_{t-k-1}+(z+w)u_{k-1}u_{t-k}\\
&=& a_1(2a_1n+a_2t+2a_3)u_{t-k}(u_k+u_{k-1})\\
& & +(b_1+b_2)(2b_1n+b_2t+2b_3)u_{k-1}(u_{t-k-1}+u_{t-k})\\
& &
+(t-2k)^2(a_2^2u_ku_{t-k}+b_2^2u_{k-1}u_{t-k-1}+2a_2b_2u_{k-1}u_{t-k}).
\end{eqnarray*}
Note that
\begin{eqnarray*}
&&2a_1n+a_2t+2a_3=[a_1n+a_2(t-k)+a_3]+(a_1n+a_2k+a_3),\\
&&2b_1n+b_2t+2b_3=[b_1n+b_2(t-k)+b_3]+(b_1n+b_2k+b_3).
\end{eqnarray*}
Hence $(2a_1n+a_2t+2a_3)u_{t-k}\ge 0$ and
$(2b_1n+b_2t+2b_3)u_{k-1}\ge 0$ for $0\le k<t-k\le n$. Moreover,
$a_2^2u_ku_{t-k}+b_2^2u_{k-1}u_{t-k-1}+2a_2b_2u_{k-1}u_{t-k}\ge 0$
by the arithmetic-geometric mean inequality and the log-concavity of
$\{u_k\}$. Consequently, we have $c_k+c_{t-k}\ge 0$ for all $0\le
t\le 2n$ and $0\le k\le\lrf{t/2}$. This gives the required result.
The proof of the theorem is complete. \end{proof}

Now the $q$-log-convexity of the Bell polynomials and the Eulerian
polynomials follows immediately from Theorem~\ref{T-qLCX}.
\begin{prop}\label{bnq}
The Bell polynomials $B_n(q)$ form a $q$-log-convex sequence.
\end{prop}
\begin{rem}\label{bn2}
An immediate consequence of Proposition~\ref{bnq} is the
log-convexity of the Bell numbers. Also, note that $2$-colored Bell
number $S_n=B_n(2)$ in Example~\ref{ex-sn}. Hence the log-convexity
of $\{S_n\}_{n\ge 0}$ follows from Proposition~\ref{bnq}.
\end{rem}
\begin{prop}\label{anq}
The Eulerian polynomials $A_n(q)$ form a $q$-log-convex sequence.
\end{prop}
\begin{rem}\label{an2}
Note that the ordered Bell number $c(n)=A_n(2)/2$ in
Example~\ref{ex-cn} by the Frobenius formula. Hence the
log-convexity of $\{c(n)\}_{n\ge 0}$ follows from
Proposition~\ref{anq}.
\end{rem}
\subsection{Linear transformations preserving log-convexity}
\hspace*{\parindent}
In \cite{WYjcta07}, Wang and Yeh established the connection between
linear transformations preserving the log-concavity and the
$q$-log-concavity. This method is also effective for the
log-convexity.

Given a triangle $\{a(n,k)\}_{0\le k\le n}$ of nonnegative real
numbers, consider the linear transformation
\begin{equation}\label{lt}
z_n=\sum_{k=0}^{n}a(n,k)x_k,\quad n=0,1,2,\ldots.
\end{equation}
For convenience, let $a(n,k)=0$ unless $0\le k\le n$. For $0\le t\le
2n$, define
\begin{equation*}
    a_k(n,t)=a(n-1,k)a(n+1,t-k)+a(n+1,k)a(n-1,t-k)-2a(n,k)a(n,t-k)
\end{equation*}
if $0\le k<t/2$, and
\begin{equation*}
    a_k(n,t)=a(n-1,k)a(n+1,k)-a^2(n,k)
\end{equation*}
if $t$ is even and $k=t/2$. Also, define
\begin{equation*}\label{q-p}
\A_n(q)=\sum_{k=0}^na(n,k)q^ k,\quad n=0,1,2,\ldots.
\end{equation*}
It is clear that if the linear transformation (\ref{lt}) preserves
the log-convexity, then for each positive number $q$, the sequence
$\{\A_n(q)\}$ is log-convex. On the other hand, we have the
following.
\begin{thm}\label{q-general}
Suppose that the triangle $\{a(n,k)\}$ of nonnegative real numbers
satisfies the following two conditions.
\begin{itemize}
  \item [\rm(C1)] The sequence of polynomials $\{\A_n(q)\}_{n\ge 0}$ is $q$-log-convex.
  \item [\rm(C2)] There exists an index $r=r(n,t)$ such that $a_k(n,t)\ge
  0$ for $k\le r$ and $a_k(n,t)<0$ for $k>r$.
\end{itemize}
Then the following two results hold.
\begin{itemize}
\item [\rm(R1)]
The linear transformation $z_n=\sum_{k=0}^{n}a(n,k)x_k$ preserves
the log-convexity.
\item [\rm(R2)]
If the sequence $\{u_k\}_{k\ge 0}$ is log-convex and
$b(n,k)=a(n,k)u_k$ for $0\le k\le n$, then the triangle
$\{b(n,k)\}_{0\le k\le n}$ also satisfies the conditions (C1) and
(C2).
\end{itemize}
\end{thm}
\begin{proof} Note that
$$z_{n-1}z_{n+1}-z_n^2=\sum_{t=0}^{2n}\left[\sum_{k=0}^{\lrf{t/2}}a_k(n,t)x_kx_{t-k}\right]$$
and
$$\A_{n-1}(q)\A_{n+1}(q)-\A_n^2(q)=\sum_{t=0}^{2n}\left[\sum_{k=0}^{\lrf{t/2}}a_k(n,t)\right]q^t.$$
Denote $\mathcal {A}(n,t)=\sum_{k=0}^{\lrf{t/2}}a_k(n,t)$. Then the
condition (C1) is equivalent to $\mathcal {A}(n,t)\ge 0$ for $0\le
t\le 2n$. Assume that $\{x_k\}_{k\ge 0}$ is log-convex. Then
$x_0x_t\ge x_1x_{t-1}\ge x_2x_{t-2}\ge\cdots$. It follows that
$$\sum_{k=0}^{\lrf{t/2}}a_k(n,t)x_kx_{t-k}\ge
\sum_{k=0}^{\lrf{t/2}}a_k(n,t)x_rx_{t-r}=\mathcal
{A}(n,t)x_rx_{t-r}$$ by the condition (C2). Thus
$z_{n-1}z_{n+1}-z_n^2\ge\sum_{t=0}^{2n}\mathcal
{A}(n,t)x_rx_{t-r}\ge 0$, and $\{z_n\}_{n\ge 0}$ is therefore
log-convex. This proves (R1).

Note that $b_k(n,t)=a_k(n,t)u_ku_{t-k}$ by the definition. Hence the
triangle $\{b(n,k)\}_{0\le k\le n}$ satisfies the condition (C2). On
the other hand, we have
$$\mathcal {B}(n,t)=\sum_{k=0}^{\lrf{t/2}}b_k(n,t)
=\sum_{k=0}^{\lrf{t/2}}a_k(n,t)u_ku_{t-k}\ge
\sum_{k=0}^{\lrf{t/2}}a_k(n,t)u_ru_{t-r}=\mathcal
{A}(n,t)u_ru_{t-r}\ge 0$$ by the log-convexity of $u_k$ and the
condition (C2). So the triangle $\{b(n,k)\}_{0\le k\le n}$ satisfies
the condition (C1). This proves (R2).\end{proof}
\begin{prop}\label{n+k--2k}
The linear transformation
$$z_n=\sum_{k=0}^{n}\binom{n+k}{n-k}x_k,\quad n=0,1,2,\ldots$$
preservers the log-convexity of sequences.
\end{prop}
\begin{proof}
Let $a(n,k)=\binom{n+k}{n-k}$ for $0\le k\le n$. Then by
Theorem~\ref{q-general}, it suffices to show that the triangle
$\{a(n,k)\}$ satisfies the conditions (C1) and (C2).

Let $\A_n(q)=\sum_{k=0}^{n}\binom{n+k}{n-k}q^k$, which is the $n$-th
Morgan-Voyce polynomial (\cite{Swa68}). By the recurrence relation
of the binomial coefficients, we can obtain
$$\binom{n+1+k}{n+1-k}=2\binom{n+k}{n-k}+\binom{n+k-1}{n-k+1}-\binom{n-1+k}{n-1-k}.$$
From this it follows that $\A_{n+1}(q)=(2+q)\A_n(q)-\A_{n-1}(q)$,
which is equivalent to
$$\begin{pmatrix}
\A_{n+1} & \A_n\\ \A_n & \A_{n-1}
\end{pmatrix}
=\begin{pmatrix} 2+q & -1\\ 1 & 0
\end{pmatrix}
\begin{pmatrix}
\A_{n} & \A_{n-1}\\ \A_{n-1} & \A_{n-2}
\end{pmatrix}.$$
Now consider the determinants on the two sides of the equality. Then
we have
$$\A_{n-1}(q)\A_{n+1}(q)-\A^2_n(q)
=\A_{n-2}(q)\A_n(q)-\A^2_{n-1}(q)
=\cdots=\A_0(q)\A_2(q)-\A^2_1(q)=q$$ by the initial conditions
$\A_0(q)=1,\A_1(q)=1+q$ and $\A_2(q)=1+3q+q^2$. Thus the sequence
$\{\A_n(q)\}_{n\ge 0}$ is $q$-log-convex, and so the condition (C1)
is satisfied.

By the definition, we have
\begin{eqnarray*}
a_k(n,t) &=& \binom{n-1+k}{n-1-k}\binom{n+1+t-k}{n+1-t+k}
+\binom{n+1+k}{n+1-k}\binom{n-1+t-k}{n-1-t+k}\\
&& -2\binom{n+k}{n-k}\binom{n+t-k}{n-t+k}\\
&=&
\frac{(n-1+k)!(n-1+t-k)!}{(2k)!(2t-2k)!(n+1-k)!(n+1-t+k)!}\widetilde{a}_k(n,t)
\end{eqnarray*}
when $k<t/2$, and
\begin{eqnarray*}
a_k(n,t) =
\binom{n-1+k}{n-1-k}\binom{n+1+k}{n+1-k}-\binom{n+k}{n-k}^2 =
\frac{1}{2}\left[\frac{(n-1+k)!}{(2k)!(n+1-k)!}\right]^2\widetilde{a}_k(n,t)
\end{eqnarray*}
when $t$ even and $k=t/2$, where
\begin{eqnarray*}
\widetilde{a}_k(n,t) &=& (n+k)(n+1+k)(n-t+k)(n-t+k+1)\\
&& +(n-k)(n-k+1)(n+t-k)(n+t-k+1)\\
&& -2(n+k)(n-k+1)(n+t-k)(n-t+k+1).
\end{eqnarray*}
Clearly, $a_k(n,t)$ has the same sign as that of
$\widetilde{a}_k(n,t)$ for each $k$. Using Maple, we obtain that the
derivative of $\widetilde{a}_k(n,t)$ with respect to $k$ is
$-2(t-2k)[2(2n+1)^2-t]\le 0$. Thus $\widetilde{a}_k(n,t)$ changes
sign at most once (from nonnegative to nonpositive), and so does
$a_k(n,t)$. Thus the condition (C2) is also satisfied. This
completes our proof.
\end{proof}

Note that the even-indexed Fibonacci numbers
$F_{2n}=\sum_{k=0}^{n}\binom{n+k}{n-k}$. Hence the log-convexity of
the numbers $F_{2n}$ follows immediately from
Proposition~\ref{n+k--2k}. It is also known that the large
Schr\"oder numbers $r_n=\sum_{k=0}^{n}\binom{n+k}{n-k}C_k$ and the
central Delannoy numbers $D(n)=\sum_{k=0}^{n}\binom{n+k}{n-k}b(k)$.
So the log-convexity of the numbers $r_n$ and $D(n)$ is implied by
the log-convexity of the Catalan numbers $C_k$ and the central
binomial coefficients $b(k)$ respectively.

The $q$-Schr\"oder number $r_n(q)$, introduced by Bonin, Shapiro and
Simion~\cite{BSS93}, is defined as the $q$-analog of the large
Schr\"oder number $r_n$:
$$r_n(q)=\sum_{P}q^{\mathrm{diag}(P)},$$
where $P$ takes over all Schr\"oder paths from $(0,0)$ to $(n,n)$
and $\mathrm{diag} (P)$ denotes the number of diagonal steps in the
path $P$. Clearly,  $r_n(1)=r_n$, the large Schr\"oder numbers.
Also,
$$r_n(q)=\sum_{k=0}^n\binom{n+k}{n-k}C_kq^{n-k}.$$
From Theorem~\ref{q-general} and Proposition~\ref{n+k--2k} we can
obtain the $q$-log-convexity of the $q$-Schr\"oder numbers.
\begin{cor}\label{rnq}
The $q$-Schr\"oder numbers $r_n(q)$ form a $q$-log-convex sequence.
\end{cor}

Similarly, consider the $q$-central Delannoy numbers
$$D_n(q)=\sum_{k=0}^{n}\binom{n+k}{n-k}b(k)q^{n-k}$$
(Sagan~\cite{Sag98}). We have the following.
\begin{cor}
The $q$-central Delannoy numbers $D_n(q)$ form a $q$-log-convex
sequence.
\end{cor}
\section{Concluding remarks and open problems}
\hspace*{\parindent}
In this paper we have explored the log-convexity of some
combinatorial sequences by algebraic and analytic approaches. It is
natural to look for combinatorial interpretations for the
log-convexity of these sequences since their strong background in
combinatorics. Callan~\cite{Cal00} gave an injective proof for the
log-convexity of the Motzkin numbers. It is possible to give
combinatorial interpretations for the log-convexity of more
combinatorial numbers. We feel that the lattice path techniques of
Wilf~\cite{Wil77} and Gessel-Viennot~\cite{GV85} are useful. As an
example, we give an injective proof for the log-convexity of the
Catalan numbers.

Recall that the Catalan number $C_n$ is the number of lattice paths
from $(i,i)$ to $(n+i,n+i)$ with steps $(0,1)$ and $(1,0)$ and never
rising above the line $y=x$ (see Stanley~\cite[Exercise 6.19
(h)]{Sta97II} for instance). Let $\mathscr{C}_n(i)$ be the set of
such paths. We next show that $C_n^2\le C_{n+1}C_{n-1}$ by
constructing an injection
$$\phi: \mathscr{C}_n(0)\times\mathscr{C}_n(1)
\rightarrow\mathscr{C}_{n+1}(0)\times\mathscr{C}_{n-1}(1).$$

Consider a path pair
$(p,q)\in\mathscr{C}_n(0)\times\mathscr{C}_n(1)$. Clearly, $p$ and
$q$ must intersect. Let $C$ be their first intersect point. Then $C$
splits $p$ into two parts $p_1$ and $p_2$, and splits $q$ into two
parts $q_1$ and $q_2$. Thus the concatenation $p'$ of $p_1$ and
$q_2$ is a path in $\mathscr{C}_{n+1}(0)$, and the concatenation
$q'$ of $q_1$ and $p_2$ is a path in $\mathscr{C}_{n-1}(1)$. Define
$\phi(p,q)=(p',q')$. Then the image of $\phi$ consists of precisely
$(p',q')\in\mathscr{C}_{n+1}(0)\times\mathscr{C}_{n-1}(1)$ such that
$p'$ and $q'$ intersect. It is easy to see that if
$\phi(p,q)=(p',q')$, then applying the same algorithm to $(p',q')$
recovers $(p,q)$. Thus $\phi$ is injective, as desired.

It would be interesting to have a combinatorial interpretation for
the log-convexity of combinatorial sequences satisfying a three-term
recurrence. We  refer the reader to Sagan~\cite{Sag88} for
combinatorial proofs for the log-concavity of combinatorial
sequences satisfying a three-term recurrence.

The log-convexity of the Bell numbers has been shown by several
different approaches in \S2. An intriguing problem is to find a
combinatorial interpretation for the log-convexity of the Bell
numbers.

The Narayana number $N(n,k)$ is defined as the number of Dyck paths
of length $2n$ with exactly $k$ peaks (a peak of a path is a place
at which the step $(1,1)$ is directly followed by the step
$(1,-1)$). The Narayana numbers have an explicit expression
$N(n,k)=\frac{1}{n}\binom{n}{k}\binom{n}{k-1}$. The Narayana
polynomials $N_n(q)=\sum_{k=0}^{n}N(n,k)q^k$ are the generating
function of the Narayana numbers and satisfy the recurrence
\begin{equation}\label{rr-nq}
(n+1)N_n(q)=(2n-1)(1+q)N_{n-1}(q)-(n-2)(1-q)^2N_{n-2}(q)
\end{equation}
(see Sulanke~\cite{Sul02} for a combinatorial proof). The Narayana
polynomials $N_n(q)$ are closely related to the $q$-Schr\"oder
numbers $r_n(q)$. It is known that $r_n(q)=N_n(1+q)$ (\cite{Sul02}).
We have showed the $q$-log-convexity of the $q$-Schr\"oder numbers
$r_n(q)$ in Corollary~\ref{rnq}. We also propose the following
stronger conjecture.
\begin{conj}
The Narayana polynomials $N_n(q)$ form a $q$-log-convex sequence.
\end{conj}

This conjecture has been verified for $n\le 100$ using Maple. It can
also be shown that for each fixed nonnegative number $q$, the
sequence $\{N_n(q)\}_{n\ge 0}$ is log-convex by means of
Theorem~\ref{c-} and the recurrence (\ref{rr-nq}). As consequences,
the Catalan numbers $C_n=\sum_{k=0}^{n}N(n,k)$ and the large
Schr\"oder number $r_n=\sum_{k=0}^{n}N(n,k)2^k$ (\cite{Sul02}) form
log-convex sequences respectively. A problem
naturally arises. 
\begin{conj}
The Narayana transformation $z_n=\sum_{k=0}^{n}N(n,k)x_k$ preserves
the log-convexity.
\end{conj}

It is also known that the central binomial coefficients
$b(n)=\sum_{k=0}^{n}\binom{n}{k}^2$ and the central Delannoy numbers
$D(n)=\sum_{k=0}^{n}\binom{n}{k}^22^k$ (\cite{Sul03}). So we propose
the following.
\begin{conj}
The triangle $\left\{\binom{n}{k}^2\right\}$ satisfies the
conditions (C1) and (C2) in Theorem~\ref{q-general}.
\end{conj}

In Proposition~\ref{anq}, we obtain the $q$-log-convexity of the
Eulerian polynomials. A closely related problem is the following.
\begin{conj}\label{Eule}
The Eulerian transformation $z_n=\sum_{k=0}^{n}A(n,k)x_k$ preserves
the log-convexity.
\end{conj}

In this paper we show that the Bell polynomials, the Eulerian
polynomials, the Morgan-Voyce polynomials, the Narayana polynomials,
the $q$-central Delannoy numbers and the $q$-Schr\"oder numbers are
$q$-log-convex respectively. All these polynomials can be shown,
using the methods established in \cite{LW-RZP,WYjcta05}, to have
only real zeros. It seems that this relation deserves further study
and investigation.
\section*{Acknowledgments}
\hspace*{\parindent}
The authors thank the anonymous referee for his/her constructive
comments and helpful suggestions which have greatly improved the
original manuscript.


\begin{thebibliography}{99}
\bibitem{Aig98}
M. Aigner, Motzkin numbers, European J. Combin. 19 (1998) 663--675.
\bibitem{BC96}
E.A. Bender, E.R. Canfield, Log-concavity and related properties of
the cycle index polynomials, J. Combin. Theory Ser. A 74 (1996)
57--70.
\bibitem{BSS93}
J. Bonin, L. Shapiro, R. Simion, Some $q$-analogues of the
Schr\"oder numbers arising from combinatorial statistics on lattice
paths, J. Statist. Plann. Inference 34 (1993) 35--55.
\bibitem{Bre89}
F. Brenti, Unimodal, log-concave, and P\'olya frequency sequences in
combinatorics,  Mem. Amer. Math. Soc. 413 (1989).
\bibitem{Bre94}
F. Brenti, Log-concave and unimodal sequences in algebra,
combinatorics, and geometry: An update, Contemp. Math., vol. 178,
1994, pp. 71--89.
\bibitem{Bre95}
F. Brenti, Combinatorics and total positivity, J. Combin. Theory
Ser. A 71 (1995) 175--218.
\bibitem{Bre96}
F. Brenti, The applications of total positivity to combinatorics,
and conversely, in: Total Positivity and Its Applications, Jaca,
1994, in: Math. Appl., vol. 359, Kluwer, Dordrecht, 1996, pp.
451--473.
\bibitem{Cal00}
D. Callan,
Notes on Motzkin and Schr\"oder numbers, \\
\verb|http://www.stat.wisc.edu/~callan/papersother/|.
\bibitem{Com74}
L. Comtet, Advanced Combinatorics, Reidel, Dordrecht, 1974.
\bibitem{DP49}
H. Davenport, G. P\'olya, On the product of two power series, Canad.
J. Math. 1 (1949) 1--5.
\bibitem{DV03}
T. Do\v sli\'c, D. Veljan, Calculus proofs of some combinatorial
inequalities, Math. Inequal. Appl. 6 (2003) 197--209.
\bibitem{DSV04}
T. Do\v sli\'c, D. Svrtan, D. Veljan, Enumerative aspects of
secondary structures, Discrete Math. 285 (2004) 67--82.
\bibitem{Dos05}
T. Do\v sli\'c, Log-balanced combinatorial sequences, Int. J. Math.
Math. Sci. 4 (2005) 507--522.
\bibitem{Dos06}
T. Do\v sli\'c, Logarithmic behavior of some combinatorial
sequences, math/CO. 0603379.
\bibitem{FZ97}
D. Foata, D. Zeilberger, A classic proof of a recurrence for a very
classical sequence, J. Combin. Theory Ser. A 80 (1997) 380--384.
\bibitem{EHPR98}
A. Ehrenfeucht, T. Harju, P. ten Pas, G. Rozenberg, Permutations,
parenthesis words, and Schr\"oder numbers, Discrete Math. 190 (1998)
259--264.
\bibitem{Eng94}
K. Engel, On the average rank of an element in a filter of the
partition lattice, J. Combin. Theory Ser. A 65 (1994) 67--78.
\bibitem{GV85}
I. Gessel, G. Viennot, Binomial determinants, path, and hook length
formulae, Adv. in Math. 58 (1985) 300--321.
\bibitem{Guy00}
R.K. Guy, Catwalks, sandsteps and Pascal pyramids, J. Integer Seq. 3
(2000) Article 00.1.6.
\bibitem{HR70}
F. Harary, R.C. Read, The enumeration of tree-like polyhexes, Proc.
Edinburgh Math. Soc. (2) 17 (1970) 1--13.
\bibitem{Kar68}
S. Karlin, Total Positivity, vol. I, Stanford University Press,
Stanford, 1968.
\bibitem{Kur72}
D.C. Kurtz, A note on concavity properties of triangular arrays of
numbers, J. Combin. Theory Ser. A 13 (1972) 135--139.
\bibitem{LW-RZP}
L.L. Liu, Y. Wang, A unified approach to polynomial sequences with
only real zeros, Adv. in Appl. Math. (2006),
doi:10.1016/j.aam.2006.02.003.
\bibitem{PW99}
P. Peart, W.-J. Woan, Bijective proofs of the Catalan and fine
recurrences, Congr. Numer. 137 (1999) 161--168.
\bibitem{PW02}
P. Peart, W.-J. Woan, A bijective proof of the Delannoy recurrence,
Congr. Numer. 158 (2002) 29--33.
\bibitem{Rom84}
S. Roman, The Umbral Calculus,  Academic Press, New York, 1984.
\bibitem{Sag88}
B.E. Sagan, Inductive and injective proofs of log concavity results,
Discrete Math. 68 (1988) 281--292.
\bibitem{Sag92TAMS}
B.E. Sagan, Log concave sequences of symmetric functions and analogs
of the Jacobi-Trudi determinants, Trans. Amer. Math. Soc. 329 (1992)
795--811.
\bibitem{Sag92DM}
B.E. Sagan, Inductive proofs of $q$-log concavity, Discrete Math. 99
(1992) 298--306.
\bibitem{Sag98}
B.E. Sagan, Unimodality and the reflection principle, Ars Combin. 48
(1998) 65--72.
\bibitem{Slo}
N.J.A. Sloane, The On-Line Encyclopedia of Integer Sequences,\\
\verb|http:// www.research.att.com/¡«njas/sequences/|.
\bibitem{Sta89}
R.P. Stanley, Log-concave and unimodal sequences in algebra,
combinatorics, and geometry, Ann. New York Acad. Sci. 576 (1989)
500--534.
\bibitem{Sta97I}
R.P. Stanley, Enumerative Combinatorics, vol. 1, Cambridge Univ.
Press, Cambridge, UK, 1997.
\bibitem{Sta97II}
R.P. Stanley, Enumerative Combinatorics, vol. 2, Cambridge Univ.
Press, Cambridge, UK, 1997.
\bibitem{Sul98}
R.A. Sulanke, Bijective recurrences concerning Schr\"oder paths,
Electron. J. Combin. 5 (1998), Research Paper 47, 11 pp.
\bibitem{Sul01}
R.A. Sulanke, Bijective recurrences for Motzkin paths, Adv. in Appl.
Math. 27 (2001) 627--640.
\bibitem{Sul02}
R.A. Sulanke, The Narayana distribution, J. Statist. Plann.
Inference 101 (2002) 311--326.
\bibitem{Sul03}
R.A. Sulanke, Objects counted by the central Delannoy numbers, J.
Integer Seq. 6 (2003) Article 03.1.5, 19 pp.
\bibitem{Swa68}
M.N.S. Swamy, Further properties of Morgan-Voyce polynomials,
Fibonacci Quart. 6 (1968) 167--175.
\bibitem{BT00}
R. Theodorescu, J.M. Borwein, Problems and Solutions: Solutions:
Moments of the Poisson distribution: 10738. Amer. Math. Monthly 107
(2000) 659.
\bibitem{wyjcta02}
Y. Wang, A simple proof of a conjecture of Simion, J. Combin. Theory
Ser. A 100 (2002) 399--402.
\bibitem{wyeujc02}
Y. Wang, Proof of a conjecture of Ehrenborg and Steingr\'imsson on
excedance statistic, European J. Combin. 23 (2002) 355--365.
\bibitem{wylaa03}
Y. Wang, Linear transformations preserving log-concavity, Linear
Algebra Appl. 359 (2003) 161--167.
\bibitem{WYeujc05}
Y. Wang, Y.-N. Yeh, Proof of a conjecture on unimodality, European
J. Combin. 26 (2005) 617--627.
\bibitem{WYjcta05}
Y. Wang, Y.-N. Yeh, Polynomials with real zeros and P\'olya
frequency sequences, J. Combin. Theory Ser. A 109 (2005) 63--74.
\bibitem{WYjcta07}
Y. Wang, Y.-N. Yeh, Log-concavity and LC-positivity, J. Combin.
Theory Ser. A 114 (2007) 195-210.
\bibitem{Wes95}
J. West, Generating trees and the Catalan and Schr\"oder numbers,
Discrete Math. 146 (1995) 247--262.
\bibitem{Wil77}
H.S. Wilf, A unified setting for sequencing, ranking, and selection
algorithms for combinatorial objects, Adv. in Math. 24 (1977)
281--291.
\end{thebibliography}
\end{document}